\documentclass[11pt,twoside]{amsart}

\usepackage{latexsym}
\usepackage{amssymb}
\usepackage{amsfonts}
\usepackage{amstext}
\usepackage{multicol}
\usepackage{mathtools}
\usepackage{amscd,amssymb,amsfonts,amsmath,amstext,amsthm,verbatim,mathdots}
\usepackage{comment,array,blkarray,multirow,tikz,graphicx,ifthen}
\usepackage [margin = 1.3in] {geometry}
\usepackage{youngtab}
\usepackage{ytableau}
\usepackage{cancel}
\usepackage{hyperref}
\numberwithin{equation}{section}

\newtheorem{thm}{Theorem}[section]
\newtheorem{lem}[thm]{Lemma}
\newtheorem{cor}[thm]{Corollary}

\theoremstyle{definition}
\newtheorem{defn}[thm]{Definition}

\newtheorem{ex}[thm]{Example}
\newtheorem{remark}[thm]{Remark}

\newcommand{\bb}[1]{\mathbb{#1}}
\newcommand{\bbZ}{\bb{Z}}
\newcommand{\bbQ}{\bb{Q}}

\newcommand{\Zn}{\bbZ/n\bbZ}
\newcommand{\mb}[1]{\mathbf{#1}}

\newcommand{\wh}[1]{\widehat{#1}}

\newcommand{\LSym}{{\rm LSym}}
\newcommand{\ds}{\displaystyle}
\newcommand{\la}{\lambda}
\newcommand{\lvar}[2]{x_{#1}^{(#2)}}

\DeclareMathOperator{\wt}{{\rm wt}}

\title{A ratio of alternants formula for loop Schur functions}
\author{Gabriel Frieden}
\date{\today}
\thanks{The author was supported in part by NSF grants DMS-1464693 and DMS-0943832.}
\address{Department of Mathematics, University of Michigan, Ann Arbor, MI 48109-1043, USA}
\email{gfrieden@umich.edu}

\begin{document}
\maketitle

\begin{abstract}
Lam and Pylyavskyy introduced loop symmetric functions as a generalization of symmetric functions. They defined loop Schur functions as generating functions over semistandard tableaux with respect to a ``colored weight,'' and they proved a Jacobi--Trudi-style determinantal formula for these generating functions. We prove that loop Schur functions can be expressed as a ratio of ``loop alternants,'' extending the analogy with Schur functions. As an application, we give a new proof of the loop version of the Murnaghan--Nakayama rule.
\end{abstract}

\ytableausetup{centertableaux}

\vspace{.05in}
\section{Introduction}
Motivated by their study of total positivity in loop groups, Lam and Pylyavskyy introduced a generalization of symmetric polynomials, which they called loop symmetric functions \cite{LPwhirl}. These are polynomials in $m$ sets of $n$ variables, which are invariant under a certain birational action of the symmetric group $S_m$. When $n=1$, the birational action reduces to the permutation of the variables $x_1, \ldots, x_m$, and loop symmetric functions are the symmetric polynomials in $m$ variables.

The ring of loop symmetric functions is generated by polynomials $e_k^{(r)}$ (resp., $h_k^{(r)}$), which are ``loop analogues'' of the elementary (resp., complete homogeneous) symmetric polynomials. In analogy with the classical setting, Lam and Pylyavskyy showed that certain determinants in the $e_k^{(r)}$ (resp., $h_k^{(r)}$) are generating functions for semistandard Young tableaux of a fixed shape, where the weight of a tableau is refined by ``coloring'' the boxes of the Young diagram with elements of $\Zn$. They defined loop Schur functions to be these determinants (or, equivalently, generating functions). Due to the importance of Schur polynomials in algebraic combinatorics, representation theory, and Schubert calculus, it is natural to look for loop generalizations of the many remarkable properties of Schur polynomials.

Schur polynomials were originally defined as the ratio of two alternating polynomials, or alternants; the Jacobi--Trudi determinantal formula and the Young tableau interpretation came later. The main result of this article is a generalization of the ratio of alternants formula for Schur polynomials to the loop setting (Theorem \ref{thm_roa}). This formula was stated without proof in \cite{ICCM}. Because the generalized alternants are defined using the birational $S_m$ action, they are rational functions, rather than polynomials as in the classical case; however, we derive an alternative ratio of alternants formula which expresses loop Schur functions as a ratio of two polynomials (Theorem \ref{thm_roa_alt}). We use the ratio of alternants formula to deduce a loop generalization of the Murnaghan--Nakayama rule (Theorem \ref{thm_mn}). This result was also stated in \cite{ICCM}, and was proved combinatorially by Ross \cite{Ross}.

Prior to the work of Lam and Pylyavskyy, the birational symmetric group action had been identified as a geometric $R$-matrix in the theory of geometric crystals \cite{Yam, KNO}. In other words, the action tropicalizes to a piecewise-linear formula for the combinatorial $R$-matrix for tensor products of affine crystals---in this case, the one-row Kirillov--Reshetikhin crystals of type $A_{n-1}^{(1)}$ \cite{HHIKTT}. This suggests that a function on tensor products of one-row crystals which is invariant under the combinatorial $R$-matrix should be the tropicalization of a ratio of loop symmetric functions. Lam and Pylyavskyy showed that the intrinsic energy, an important function in affine crystal theory, is in fact the tropicalization of a certain loop Schur function \cite{LPenergy} (this loop Schur function turns out to be related to our alternant formulas; see Remark \ref{rem_Jucis_formula}). Additionally, tensor products of one-row crystals can be viewed as states in the (generalized) Box-Ball system, a well-studied cellular automaton that exhibits soliton behavior \cite{HHIKTT}. Formulas for the scattering of a given state into solitons are conjecturally given by tropicalizations of a cylindric variant of loop Schur functions \cite{LPS}.

We note that the birational $S_m$ action also arose in the context of the local Langlands program \cite{BravKazh}, and was studied in \cite{Et}. Loop symmetric functions have also found application in Gromov--Witten/Donaldson--Thomas theory \cite{RossZong}; this was Ross' motivation for proving the loop Murnaghan--Nakayama rule \cite{Ross}.

This paper is organized as follows. In Section \ref{sec_background}, we review the basics of loop symmetric functions and the birational $S_m$ action. Section \ref{sec_roa} contains statements and proofs of the two ratio of alternants formulas (Theorems \ref{thm_roa} and \ref{thm_roa_alt}), and Section \ref{sec_mn} discusses the loop Murnaghan--Nakayama rule (Theorem \ref{thm_mn}).

\subsection*{Acknowledgments}

I would like to thank my advisor, Thomas Lam, for introducing me to loop symmetric functions and their connection to affine crystals.


\section{Loop symmetric functions}
\label{sec_background}

\subsection{Loop elementary and homogeneous symmetric functions}
\label{sec_loop_sym_func}

Fix integers $m,n \geq 1$. For $i \in \{1, \ldots, m\}$ and $j \in \Zn$, let $x_i^{(j)}$ be an indeterminate, and let $\bb{Q}(x_i^{(j)})$ be the field of rational functions in these $mn$ indeterminates. Let $\mb{x}_i = (x_i^{(1)}, \ldots, x_i^{(n)})$. We view the superscript as a ``color.''

Lam and Pylyavskyy \cite{LPwhirl} introduced \emph{loop symmetric functions} as a class of polynomials in the $m$ ``variables'' $\mb{x}_1, \ldots, \mb{x}_m$. The sense in which these polynomials are ``symmetric'' is discussed in \S \ref{sec_birat_action}. Here we recall the definitions of several types of loop symmetric functions.
\begin{defn}
\label{defn_e_h}
For $k \geq 1$ and $r \in \Zn$, define the \emph{loop elementary symmetric function}
\[
e_k^{(r)} = e_k^{(r)}(\mb{x}_1, \ldots, \mb{x}_m) = \sum_{1 \leq i_1 < \ldots < i_k \leq m} x_{i_1}^{(r)} x_{i_2}^{(r+1)} \cdots x_{i_k}^{(r+k-1)},
\]
and define the \emph{loop homogeneous symmetric function}
\[
h_k^{(r)} = h_k^{(r)}(\mb{x}_1, \ldots, \mb{x}_m) = \sum_{1 \leq i_1 \leq \ldots \leq i_k \leq m} x_{i_1}^{(r)} x_{i_2}^{(r-1)} \cdots x_{i_k}^{(r-k+1)}.
\]
Also set $e_0^{(r)} = h_0^{(r)} = 1$, and $e_k^{(r)} = h_k^{(r)} = 0$ for $k < 0$.
\end{defn}

Note that if we ignore colors (or take $n = 1$), then the loop elementary (resp., homogeneous) symmetric functions are simply the ordinary elementary (resp., homogeneous) symmetric polynomials in $m$ variables.

Define the \emph{ring of loop symmetric functions} (in $m$ variables), denoted $\LSym_m$, to be the subring of $\mathbb{Q}(x_i^{(j)})$ generated by the loop elementary symmetric functions. It follows from Theorem \ref{thm_JT} below that $\LSym_m$ is also generated by the loop homogeneous symmetric functions.

\subsection{Loop Schur functions}

Given partitions $\la$ and $\mu$, write $\mu \subset \la$ if the Young diagram of $\mu$ is contained in that of $\la$. If $\mu \subset \la$, define $\la/\mu$ to be the \emph{skew diagram} (or \emph{skew shape}) obtained by removing the boxes of $\mu$ from the Young diagram of $\la$. A \emph{semistandard tableau of shape $\la/\mu$} is a filling of the skew diagram $\la/\mu$ with positive integers so that each row is weakly increasing, and each column is strictly increasing. For a box $(i,j)$ in $\la/\mu$, define its \emph{content} to be $c(i,j) = i-j$ (mod $n$). For $r \in \Zn$ and a semistandard tableau $T$ of shape $\la/\mu$, define the \emph{$r$-weight} of $T$ by
\[
\wt_r(T) = \prod_{(i,j) \in \la/\mu} x_{T(i,j)}^{(r+c(i,j))}.
\]

\begin{defn}
\label{defn_Schur}
For $r \in \Zn$ and partitions $\mu \subset \la$, define the \emph{loop skew Schur function}
\[
s_{\la/\mu}^{(r)} = s_{\la/\mu}^{(r)}(\mb{x}_1, \ldots, \mb{x}_m) = \sum_T \wt_r(T)
\]
where the sum is over semistandard tableaux of shape $\la/\mu$ with entries in $\{1, \ldots, m\}$. If $\mu = \emptyset$, then $s_\la^{(r)} = s_{\la/\emptyset}^{(r)}$ is a \emph{loop Schur function}. Set $s_\emptyset^{(r)} = 1$. Note that $s_{\la/\mu}^{(r)} = 0$ if any column of $\la/\mu$ has more than $m$ boxes.
\end{defn}

\begin{ex}
\label{ex_Schur}
Let $n = 3, m=2,$ and $\lambda = (3,2)$. There are two semistandard tableaux of shape $\lambda$ with entries in $\{1,2\}$:
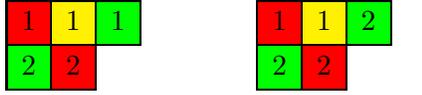
\begin{equation*}
\ytableaushort{{*(red)1}{*(yellow)1}{*(green)1},{*(green)2}{*(red)2}}
\quad\quad\quad\quad
\ytableaushort{{*(red)1}{*(yellow)1}{*(green)2},{*(green)2}{*(red)2}} \quad .
\end{equation*}
Computing 2-weights, we have
\[
s_{(3,2)}^{(2)} = \lvar{1}{2}\lvar{1}{1}\lvar{1}{3}\lvar{2}{3}\lvar{2}{2} + \lvar{1}{2}\lvar{1}{1}\lvar{2}{3}\lvar{2}{3}\lvar{2}{2}.
\]
\end{ex}

Observe that $s_{(k)}^{(r)} = h_k^{(r)}$ and $s_{(1^k)}^{(r)} = e_k^{(r)}$. These identities are special cases of Jacobi--Trudi formulas for loop skew Schur functions.

\begin{thm}[{\cite[\S 7.2]{LPwhirl}}]
\label{thm_JT}
For two partitions $\mu \subset \la$, we have
\[
s_{\la/\mu}^{(r)} = \det(h_{\la_i -\mu_j - i + j}^{(r -\mu_j + j -1)}) = \det(e^{(r+\mu'_j -j +1)}_{\la'_i - \mu'_j - i + j})
\]
where $\la'$ is the transpose of $\la$. Thus, $s_{\la/\mu}^{(r)}$ is a loop symmetric function.
\end{thm}

\begin{remark}
\label{rem_notation}
Two rings of loop symmetric functions (the ``whirl'' ring and the ``curl'' ring) are defined in \cite{LPwhirl}. We have chosen to use the ``whirl'' version, as in \cite{ICCM, LPenergy, LPS}. The ``curl'' loop Schur functions are defined by assigning to a box $(i,j)$ the content $j-i$ rather than $i-j$. In particular, what we denote by $h_k^{(r)}(\mb{x}_1, \ldots, \mb{x}_m)$ is $h_k^{(r-k+1)}(\mb{y}_1, \ldots, \mb{y}_m)$, where $\mb{y}_i = \mb{x}_{m-i+1}$, in the notation of \cite{LPwhirl}. Lam and Pylyavskyy's proof of Theorem \ref{thm_JT} uses the ``curl'' ring, so for the reader's convenience, we present the proof in ``whirl'' notation.
\end{remark}

\begin{proof}[Prof of Theorem \ref{thm_JT}]
Both identities are proved by the Lindstr\"om/Gessel--Viennot lattice path method (we refer the reader who is unfamiliar with this method to \cite[\S 2.7]{EC1}). To evaluate the first determinant, consider the network on the vertex set $\mathbb{Z}^2$, with horizontal edges from $(a-1,b)$ to $(a,b)$ of weight $x_b^{(r-a)}$, and vertical edges from $(a,b)$ to $(a,b+1)$ of weight 1. The monomials in $h_{\la_i - \mu_j - i + j}^{(r-\mu_j + j -1)}$ are the weights of the (directed) paths from $(\mu_j-j,1)$ to $(\la_i-i,m)$, so by Lindstr\"om/Gessel--Viennot, the determinant is equal to the sum over vertex-disjoint collections of $m$ paths, where the $i^{th}$ path goes from $(\mu_i - i, 1)$ to $(\lambda_i - i, m)$. Such collections of paths are in bijection with semistandard tableaux of shape $\la/\mu$ with entries in $\{1, \ldots, m\}$, and by our choice of edge weights, the weight of a collection of paths is the $r$-weight of the corresponding tableau. To evaluate the other determinant, replace the horizontal edges with diagonal edges, and argue similarly.
\end{proof}

Unfortunately, loop Schur functions do not span the ring of loop symmetric functions, and they are not even linearly independent! It remains an open problem to find a ``Schur-like'' basis of $\LSym_m$.

\subsection{Whirls and curls}
\label{sec_whirls}

Elementary (resp., homogeneous) symmetric polynomials are the coefficients of powers of $t$ in the power series $\prod_i (1+x_i t)$ (resp., $\prod_i (1-x_it)^{-1}$). Using the correspondence between power series and infinite Toeplitz matrices, these polynomials can be viewed as matrix entries of an infinite upper triangular matrix which is constant along each diagonal. Here we give an analogous description of loop elementary and homogeneous symmetric functions, following \cite{LPwhirl}.

An {\em $n$-periodic matrix} is a $\mathbb{Z} \times \mathbb{Z}$ array $(X_{ij})_{(i,j) \in \mathbb{Z}}$ such that $X_{ij} = 0$ if $i-j$ is sufficiently large, and $X_{ij} = X_{i+n,j+n}$ for all $i,j$. Multiplication of these matrices is defined in the usual way: if $X = (X_{ij})$ and $Y = (Y_{ij})$, then
\[
(XY)_{ij} = \sum_{k \in \mathbb{Z}} X_{ik}Y_{kj}.
\]
The hypothesis that $X_{ij} = 0$ for $i-j$ sufficiently large ensures that each of these sums is finite, so the product is well-defined. It's clear that the product of two $n$-periodic matrices is $n$-periodic.

Given an $n$-periodic matrix $X = (X_{ij})$, define $X^c = (-1)^{i+j} X_{ij}$. If $X$ is invertible, define $X^{-c} = (X^{-1})^c$. It is easy to see that $(XY)^c = X^cY^c$, and $(X^c)^{-1} = X^{-c}$.

\begin{defn}
Let $a^{(1)}, ..., a^{(n)}$ be indeterminates. The \emph{whirl} $M(a^{(1)}, ..., a^{(n)})$ is the $n$-periodic matrix $(M_{ij})$ with $M_{ii} = 1$, $M_{i,i+1} = a^{(i)}$ (interpret the superscripts mod $n$), and all other entries zero. The \emph{curl} $N(a^{(1)}, ..., a^{(n)})$ is the $n$-periodic matrix $N(a^{(1)}, ..., a^{(n)}) = M(a^{(1)}, ..., a^{(n)})^{-c}$.
\end{defn}

For example, when $n=3$,
\begin{equation*}
M(a^{(1)}, a^{(2)}, a^{(3)}) = \left(
\begin{array}{ccccccccc}
1 & a^{(1)} & 0 & 0 & 0 \\
0 & 1 & a^{(2)} & 0 & 0  \\
0 & 0 & 1 & a^{(3)} & 0 & \hdots \\
0 & 0 & 0 & 1 & a^{(1)}  \\
0 & 0 & 0 & 0 & 1 \\
&&\vdots &&& \ddots
\end{array}
\right)
\end{equation*}
and
\begin{equation*}
N(a^{(1)}, a^{(2)}, a^{(3)}) = \left(
\begin{array}{ccccccccc}
1 & a^{(1)} & a^{(2)}a^{(1)} & a^{(3)}a^{(2)}a^{(1)} & a^{(1)}a^{(3)}a^{(2)}a^{(1)} \\
0 & 1 & a^{(2)} & a^{(3)}a^{(2)} & a^{(1)}a^{(3)}a^{(2)}  \\
0 & 0 & 1 & a^{(3)} & a^{(1)}a^{(3)} & \hdots \\
0 & 0 & 0 & 1 & a^{(1)}  \\
0 & 0 & 0 & 0 & 1 \\
&&\vdots &&& \ddots
\end{array}
\right).
\end{equation*}
Note that we are depicting only the quadrant of the matrix with $i,j \geq 1$.

It is straightforward to show that loop elementary (resp., homogeneous) symmetric functions are the matrix entries of a product of whirls (resp., curls).

\begin{lem}[{\cite[\S 7.2]{LPwhirl}}]
\label{lem_matrix_entries}
Set
\[
A = M(\mathbf{x}_1) M(\mathbf{x}_2) \cdots M(\mathbf{x}_m), \quad\quad B = A^{-c} = N(\mathbf{x}_m) N(\mathbf{x}_{m-1}) \cdots N(\mathbf{x}_1).
\]
Then
\[
A_{ij} = e_{j-i}^{(i)}(\mb{x}_1, \ldots, \mb{x}_m) \quad\quad \text{ and } \quad\quad B_{ij} = h_{j-i}^{(j-1)}(\mb{x}_1, \ldots, \mb{x}_m).\footnote{The discrepancy between this formula and that of \cite[Lem. 7.3]{LPwhirl} is explained by Remark \ref{rem_notation}.}
\]
\end{lem}


\begin{ex}
\label{ex_loop_elem}
When $n=2$ and $m=3$,
\[
M(\mathbf{x}_1)M(\mathbf{x}_2)M(\mathbf{x}_3) =
\]
\begin{equation*}
\left(
\begin{array}{ccccccccccccccccc}
1 & \lvar{1}{1} + \lvar{2}{1} + \lvar{3}{1} & \lvar{1}{1}\lvar{2}{2} + \lvar{1}{1}\lvar{3}{2} + \lvar{2}{1}\lvar{3}{2} & \lvar{1}{1}\lvar{2}{2}\lvar{3}{1} & \\
0 & 1 & \lvar{1}{2} + \lvar{2}{2} + \lvar{3}{2} & \lvar{1}{2}\lvar{2}{1} + \lvar{1}{2}\lvar{3}{1} + \lvar{2}{2}\lvar{3}{1} & \hdots \\
0 & 0 & 1 & \lvar{1}{1} + \lvar{2}{1} + \lvar{3}{1} &  \\
0 & 0 & 0 & 1 &  \\
 &  &  \vdots &  & \ddots
\end{array}
\right).
\end{equation*}
The entries in the top row of this matrix are $e_0^{(1)}, e_1^{(1)}, e_2^{(1)}, e_3^{(1)}$ (followed by $e_4^{(1)} = e_5^{(1)} = \cdots = 0$), in agreement with Lemma \ref{lem_matrix_entries}.
\end{ex}

\begin{remark}
Theorem \ref{thm_JT} (combined with Lemma \ref{lem_matrix_entries}) shows that loop skew Schur functions are precisely the minors of the matrix $N(\mb{x}_m) \cdots N(\mb{x}_1)$. Loop Schur functions are the minors using consecutive columns; minors using consecutive rows are the loop skew Schur functions of anti-partition shape (i.e., of shape $\la/\mu$ where $\la$ is a rectangle). A similar statement can be made about minors of $M(\mb{x}_1) \cdots M(\mb{x}_m)$.\footnote{One must reflect a submatrix of $M(\mb{x}_1) \cdots M(\mb{x}_m)$ over the anti-diagonal (or equivalently, transpose and rotate the submatrix $180^\circ$) to obtain the matrix of loop elementary symmetric functions appearing in Theorem \ref{thm_JT}. Thus, loop Schur functions are the minors of $M(\mb{x}_1) \cdots M(\mb{x}_m)$ using consecutive rows.}
\end{remark}

\subsection{The birational $S_m$ action}
\label{sec_birat_action}

Symmetric polynomials in $m$ variables are the invariants of the natural action of $S_m$ on the polynomial ring in $m$ variables. Loop symmetric functions are the invariants of a more complicated $S_m$ action on $\bbQ(x_i^{(j)})$, as we now explain.

Let $\mathbf{x} = (x^{(1)}, \ldots, x^{(n)})$ and $\mathbf{y} = (y^{(1)}, \ldots, y^{(n)})$. For $r \in \Zn$, set
\[
\kappa^{(r)}(\mb{x},\mb{y}) = h_{n-1}^{(r-1)}(\mb{x}, \mb{y}) = \sum_{s = 0}^{n-1} x^{(r-1)} \cdots x^{(r-s)} y^{(r-s+1)} \cdots y^{(r-n+1)}
\]
where, as usual, the superscripts live in $\Zn$.

\begin{defn}
\label{defn_birat_action}
For $i = 1, \ldots , m-1$, define $s_i : \mathbb{Q}(x_i^{(j)}) \rightarrow \mathbb{Q}(x_i^{(j)})$ to be the $\mathbb{Q}$-algebra homomorphism which fixes $x_k^{(j)}$ for $k \neq i,i+1$, and acts on $x_i^{(j)}, x_{i+1}^{(j)}$ by
\[
s_i (x_i^{(j)}) = x_{i+1}^{(j+1)} \frac{\kappa^{(j+1)}(\mb{x}_i, \mb{x}_{i+1})}{\kappa^{(j)}(\mb{x}_i, \mb{x}_{i+1})}, \qquad s_i (x_{i+1}^{(j)}) = x_{i}^{(j-1)} \frac{\kappa^{(j-1)}(\mb{x}_i, \mb{x}_{i+1})}{\kappa^{(j)}(\mb{x}_i, \mb{x}_{i+1})}.
\]
\end{defn}

\begin{thm}[{\cite[\S 2]{Yam} \cite[\S 6]{LPwhirl}}]
\label{thm_birat_action}
\
\begin{enumerate}
\item For $i = 1, \ldots, m-1$,
\begin{equation*}
M(\mathbf{x}_i)M(\mathbf{x}_{i+1}) = s_i(M(\mathbf{x}_i)M(\mathbf{x}_{i+1})).
\end{equation*}
\item The maps $s_i$ satisfy the relations of the adjacent transpositions in the symmetric group $S_m$.
\end{enumerate}
\end{thm}

\begin{remark}
The deepest part of Theorem \ref{thm_birat_action} is the fact that the maps $s_i$ and $s_{i+1}$ satisfy the braid (or Yang--Baxter) relation. Many different proofs of this result have appeared in the literature: in addition to \cite{Yam, LPwhirl}, see \cite{Et}, \cite[\S 6.5]{LPnet}, \cite[\S 5.2]{FriR}, and the combination of \cite[\S 8.7]{BravKazh} and \cite[\S 6.2]{BerKazh}.
\end{remark}

Theorem \ref{thm_birat_action}(2) shows that the maps $s_i$ generate an action of $S_m$ on $\bbQ(x_i^{(j)})$. We call this the \emph{birational $S_m$ action} (it is also called the birational $R$-matrix in the literature). Theorem \ref{thm_birat_action}(1) and Lemma \ref{lem_matrix_entries} imply that the loop elementary symmetric functions (and thus all loop symmetric functions) are invariant under this action. In fact, the loop elementary symmetric functions are algebraically independent generators of the subring of polynomial invariants for this action \cite{LPplus}, but we will not use this result.

\subsection{Loop power sums}
\label{sec_power_sums}

We will need one additional class of loop symmetric functions. For $i = 1, \ldots, m$, set
\[
\pi_i = x_i^{(1)} x_i^{(2)} \cdots x_i^{(n)}.
\]

\begin{lem}
\label{lem_pi}
For $w \in S_m$, we have $w(\pi_i) = \pi_{w(i)}$.
\end{lem}

\begin{proof}
Since $S_m$ is generated by the $s_i$, it suffices to show that $s_i(\pi_j) = \pi_{s_i(j)}$ for each $i, j$. We compute
\begin{align*}
s_i(\pi_i) = s_i(x_i^{(1)}) s_i(x_i^{(2)}) \cdots s_i(x_i^{(n)}) = x_{i+1}^{(2)} \frac{\kappa_i^{(2)}}{\kappa_i^{(1)}}x_{i+1}^{(3)} \frac{\kappa_i^{(3)}}{\kappa_i^{(2)}} \cdots x_{i+1}^{(1)} \frac{\kappa_i^{(1)}}{\kappa_i^{(n)}} = \pi_{i+1}.
\end{align*}
Similarly, $s_i(\pi_{i+1}) = \pi_i$, and clearly $s_i(\pi_j) = \pi_j$ if $j \neq i, i+1$.
\end{proof}

\begin{defn}
\label{defn_power_sum}
For each positive integer $k$, define the \emph{loop power sum symmetric function}
\[
p_k = p_k(\mathbf{x}_1, \ldots, \mathbf{x}_m) = \sum_{i=1}^m \pi_i^k.
\]
\end{defn}

By Lemma \ref{lem_pi}, the polynomials $p_k$ are invariant under the birational $S_m$ action. The loop Murnaghan--Nakayama rule (Theorem \ref{thm_mn}) expresses $p_k$ as an alternating sum of loop Schur functions.


\section{Alternants}
\label{sec_roa}

Schur polynomials were originally defined by the formula $s_\lambda(x_1, \ldots, x_m) = a_{\lambda + \delta}/a_\delta$, where $a_\alpha$ is the determinant of the matrix $(x_j^{\alpha_i})$, and $\delta$ is the staircase partition $(m-1, m-2, \ldots, 1, 0)$. These determinants are called \emph{alternants} because they are anti-symmetric (or alternating) with respect to permutation of the variables. We now present a generalization of alternants to the loop setting.

For a sequence $\alpha = (\alpha_1, ..., \alpha_m)$ of non-negative integers, define an $m \times m$ matrix $A_{\alpha}^{(r)}$ by
\[
(A_{\alpha}^{(r)})_{ij} = t_{j, m}(x_m^{(r+m-1)}x_m^{(r+m-2)} \cdots x_m^{(r+m- \alpha_i)})
\]
where $t_{a,b}$ is the transposition in $S_m$ which swaps $a$ and $b$, acting on $\mathbb{Q}(x_i^{(j)})$ by the birational action. For example, if $n=3$ and $m=2$, then
\begin{equation}
\label{eq_A_example}
A^{(2)}_{(4,2)} =
\left(
\begin{array}{ccccccc}
\lvar{1}{2}\lvar{1}{1}\lvar{1}{3}\lvar{1}{2} \dfrac{\lvar{1}{1}\lvar{1}{3} + \lvar{1}{1}\lvar{2}{3} + \lvar{2}{1}\lvar{2}{3}}{\lvar{1}{2}\lvar{1}{1} + \lvar{1}{2}\lvar{2}{1} + \lvar{2}{2}\lvar{2}{1}} & \lvar{2}{3}\lvar{2}{2}\lvar{2}{1}\lvar{2}{3} \\
&\\
\lvar{1}{2}\lvar{1}{1} \dfrac{\lvar{1}{3}\lvar{1}{2} + \lvar{1}{3}\lvar{2}{2} + \lvar{2}{3}\lvar{2}{2}}{\lvar{1}{2}\lvar{1}{1} + \lvar{1}{2}\lvar{2}{1} + \lvar{2}{2}\lvar{2}{1}} & \lvar{2}{3}\lvar{2}{2}
\end{array}
\right).
\end{equation}
Set $a_{\alpha}^{(r)} = \det(A_{\alpha}^{(r)})$. It is easy to see that $a_\alpha^{(r)}$ is anti-symmetric with respect to the birational $S_m$ action; we call this determinant a \emph{loop alternant}. The following result was stated without proof in \cite{ICCM}.\footnote{The original statement of this result (\cite[Thm. 5.6]{ICCM}) is incorrect; using the indexing conventions of that paper, the superscript of the loop Schur function should be $r-m+1$ rather than $r-1$.}

\begin{thm}
\label{thm_roa}
For $\la$ a partition with at most $m$ parts, we have
\begin{equation}
\label{eq_roa}
s_{\lambda}^{(r)}(\mathbf{x}_1, \ldots, \mathbf{x}_m) = \frac{a_{\lambda + \delta}^{(r)}}{a_{\delta}^{(r)}},
\end{equation}
where $\delta = (m-1, m-2, \ldots, 1, 0)$.
\end{thm}

Before proving this result, we derive several corollaries. Recall that $\pi_i = x_i^{(1)} \cdots x_i^{(n)}$.

\begin{cor}
\label{cor_denom}
For $r \in \Zn$, we have
\begin{equation}
\label{eq_denom}
a_\delta^{(r)} = \dfrac{\ds \prod_{1 \leq i < j \leq m} (\pi_i-\pi_j)}{s^{(r)}_{(n-1)\delta}(\mb{x}_1, \ldots, \mb{x}_m)}
\end{equation}
where $k\delta = (k(m-1), k(m-2), \ldots, k, 0)$.
\end{cor}

\begin{proof}
Taking $\lambda = (n-1)\delta$ in Theorem \ref{thm_roa}, we have
\begin{equation*}
a_{\delta}^{(r)} = \dfrac{a_{n\delta}^{(r)}}{s_{(n-1) \delta}^{(r)}(\mathbf{x}_1, \ldots, \mathbf{x}_m)}.
\end{equation*}
The last column of $A_{n\delta}^{(r)}$ has entries $(A_{n\delta}^{(r)})_{i,m} = \pi_m^{m-i}$. Since the birational $S_m$ action permutes the $\pi_j$, we have $(A_{n\delta}^{(r)})_{i,j} = \pi_j^{m-i}$, so $A_{n\delta}^{(r)}$ has determinant $\prod_{i < j} (\pi_i - \pi_j)$.
\end{proof}

We now combine Theorem \ref{thm_roa} and Corollary \ref{cor_denom} with results of Lam and Pylyavskyy to obtain a variant of Theorem \ref{thm_roa} which expresses $s_\la^{(r)}$ as a ratio of two polynomials.

For $k \geq 1$ and $r \in \Zn$, let
\[
\sigma_k^{(r)}(\mb{x}_a, \ldots, \mb{x}_b) = \sum x_{i_1}^{(r)} x_{i_2}^{(r-1)} \cdots x_{i_k}^{(r-k+1)},
\]
where the sum is over weakly increasing sequences $a \leq i_1 \leq \cdots \leq i_k \leq b$, such that each of the numbers $a+1, a+2, \ldots, b$ appears in the sequence at most $n-1$ times. Set $\sigma_0^{(r)} = 1$. For example, if $n = 3$, then
\[
\sigma_4^{(2)}(\mb{x}_1, \mb{x}_2) = \lvar{1}{2} \lvar{1}{1}\lvar{1}{3} \lvar{1}{2} + \lvar{1}{2} \lvar{1}{1}\lvar{1}{3} \lvar{2}{2} + \lvar{1}{2} \lvar{1}{1}\lvar{2}{3} \lvar{2}{2} = \lvar{1}{2} \lvar{1}{1} \sigma_2^{(3)}(\mb{x}_1, \mb{x}_2).
\]

For a sequence $\alpha = (\alpha_1, \ldots, \alpha_m)$ of non-negative integers, define an $m \times m$ matrix $B_\alpha^{(r)}$ by
\[
(B_\alpha^{(r)})_{ij} = x_j^{(r+j-1)} x_j^{(r+j-2)} \cdots x_j^{(r+j-\alpha_i)} \sigma_{(n-1)(m-j)}^{(r+j-\alpha_i-1)}(\mb{x}_j, \ldots, \mb{x}_m).
\]
Set $b_\alpha^{(r)} = \det(B_{\alpha}^{(r)})$.

\begin{thm}
\label{thm_roa_alt}
For $\la$ a partition with at most $m$ parts, we have
\begin{equation}
\label{eq_roa_alt}
s_\la^{(r)}(\mb{x}_1, \ldots, \mb{x}_m) = \dfrac{b_{\la + \delta}^{(r)}}{b_{\delta}^{(r)}} = \dfrac{b_{\la + \delta}^{(r)}}{\ds \prod_{1 \leq i < j \leq m} (\pi_i - \pi_j)}.
\end{equation}
\end{thm}

\begin{ex}
\label{ex_B}
Let $n=3, m = 2,$ and $\lambda = (3,2)$. Then
\begin{equation*}
B^{(2)}_{\lambda + \delta} = B^{(2)}_{(4,2)} =
\left(
\begin{array}{ccccccc}
\lvar{1}{2}\lvar{1}{1}\lvar{1}{3}\lvar{1}{2} (\lvar{1}{1}\lvar{1}{3} + \lvar{1}{1}\lvar{2}{3} + \lvar{2}{1}\lvar{2}{3}) & \lvar{2}{3}\lvar{2}{2}\lvar{2}{1}\lvar{2}{3} \\
&\\
\lvar{1}{2}\lvar{1}{1} (\lvar{1}{3}\lvar{1}{2} + \lvar{1}{3}\lvar{2}{2} + \lvar{2}{3}\lvar{2}{2}) & \lvar{2}{3}\lvar{2}{2}
\end{array}
\right)
\end{equation*}
(cf. \eqref{eq_A_example}). The determinant of this matrix is
\[
b^{(2)}_{(4,2)} = \lvar{1}{2}\lvar{1}{1}\lvar{2}{3}\lvar{2}{2}(\lvar{1}{3}\lvar{1}{2}\lvar{1}{1}\lvar{1}{3} + \lvar{1}{3}\lvar{1}{2}\lvar{1}{1}\lvar{2}{3} - \lvar{1}{3}\lvar{2}{2}\lvar{2}{1}\lvar{2}{3} - \lvar{2}{3}\lvar{2}{2}\lvar{2}{1}\lvar{2}{3}).
\]
Similarly, one computes $b^{(2)}_{(1,0)} = x_1^{(1)}x_1^{(2)}x_1^{(3)} - x_2^{(1)}x_2^{(2)}x_3^{(3)}$, and the ratio $b^{(2)}_{(4,2)}/b^{(2)}_{(1,0)}$ is indeed equal to the loop Schur function $s^{(2)}_{(3,2)}$ from Example \ref{ex_Schur}.
\end{ex}

\begin{proof}[Proof of Theorem \ref{thm_roa_alt}]
The entries of the matrix $A_\alpha^{(r)}$ are given by
\begin{align*}
(A_{\alpha}^{(r)})_{ij} &= t_{j, m}(x_m^{(r+m-1)}x_m^{(r+m-2)} \cdots x_m^{(r+m- \alpha_i)}) \\
&= s_j s_{j+1} \cdots s_{m-2} s_{m-1} (x_m^{(r+m-1)} x_m^{(r+m-2)} \cdots x_m^{(r+m- \alpha_i)}),
\end{align*}
where the second equality comes from writing $t_{j,m} = s_j s_{j+1} \cdots s_{m-2} s_{m-1} s_{m-2} \cdots s_{j+1} s_j$, and observing that the maps $s_j, \ldots, s_{m-2}$ do not affect the variables $x_m^{(a)}$. By \cite[Lem. 3.1]{LPenergy},
\[
s_j s_{j+1} \cdots s_{m-2}s_{m-1}(x_m^{(a)}) = x_j^{(a-m+j)} \dfrac{\sigma_{(n-1)(m-j)}^{(a-m+j-1)}(\mb{x}_j, \ldots, \mb{x}_m)}{\sigma_{(n-1)(m-j)}^{(a-m+j)}(\mb{x}_j, \ldots, \mb{x}_m)},
\]
so
\[
(A_{\alpha}^{(r)})_{ij} = x_j^{(r+j-1)} x_j^{(r+j-2)} \cdots x_j^{(r+j-\alpha_i)} \dfrac{\sigma_{(n-1)(m-j)}^{(r+j-\alpha_i-1)}(\mb{x}_j, \ldots, \mb{x}_m)}{\sigma_{(n-1)(m-j)}^{(r+j-1)}(\mb{x}_j, \ldots, \mb{x}_m)} = \dfrac{(B_\alpha^{(r)})_{ij}}{\sigma_{(n-1)(m-j)}^{(r+j-1)}(\mb{x}_j, \ldots, \mb{x}_m)}.
\]
Taking determinants, we have
\begin{equation}
\label{eq_a_b}
a_{\alpha}^{(r)} = \dfrac{b_\alpha^{(r)}}{\prod_{j=1}^m \sigma_{(n-1)(m-j)}^{(r+j-1)}(\mb{x}_j, \ldots, \mb{x}_m)}.
\end{equation}
This identity holds for all $\alpha \in (\bbZ_{\geq 0})^m$, so
\[
\dfrac{b^{(r)}_{\lambda+\delta}}{b^{(r)}_\delta} = \dfrac{a^{(r)}_{\lambda + \delta}}{a^{(r)}_\delta} = s^{(r)}_\lambda(\mb{x}_1, \ldots, \mb{x}_m)
\]
by Theorem \ref{thm_roa}.

By \cite[Thm. 2.5]{LPenergy}, the denominator of \eqref{eq_a_b} is equal to $s_{(n-1)\delta}^{(r)}(\mb{x}_1, \ldots, \mb{x}_m)$. The second equality in \eqref{eq_roa_alt} now follows from Corollary \ref{cor_denom}.
\end{proof}

\begin{proof}[Proof of Theorem \ref{thm_roa}]
This proof is adapted from an argument in \cite[\S I.3]{Mac}. Unless otherwise noted, all loop symmetric functions are in the variables $\mathbf{x}_1, \ldots, \mathbf{x}_m$.

Set $\alpha = \lambda + \delta$. Define an $m \times m$ matrix $H_{\alpha}^{(r)}$ by
\[
(H_{\alpha}^{(r)})_{ij} = h_{\alpha_i - m + j}^{(r+j-1)} = h_{\lambda_i - i + j}^{(r+j-1)}.
\]
Note that $\det(H_{\alpha}^{(r)}) = s_{\lambda}^{(r)}$ by Theorem \ref{thm_JT}. Define an $m \times m$ matrix $E^{(r)}$ by
\[
(E^{(r)})_{ij} = (-1)^{m-i}t_{j,m}(e_{m-i}^{(r+i)}[\wh{m}])
\]
where $e_k^{(r)}[\wh{m}] = e_k^{(r)}(\mathbf{x}_1, \ldots, \mathbf{x}_{m-1})$. We will show that
\begin{equation}
\label{eq_ultimatewant}
H_{\alpha}^{(r)} E^{(r)} = A_{\alpha}^{(r)}.
\end{equation}
To this end, set $X = N(\mathbf{x}_m) \cdots N(\mathbf{x}_1)$ and $Y = M(\mathbf{x}_1)^c \cdots M(\mathbf{x}_{m-1})^c$. Clearly $XY = N(\mathbf{x}_m)$, so
\begin{equation}
\label{eq_ab}
(XY)_{uv} = \sum_k X_{uk}Y_{kv} = h_{v-u}^{(v-1)}(\mathbf{x}_m)
\end{equation}
by Lemma \ref{lem_matrix_entries}. Applying $t_{j,m}$ to \eqref{eq_ab} and using the fact that $t_{j,m}$ fixes every matrix entry of $X$, we obtain
\begin{equation}
\label{eq_t_jm}
\sum_k X_{uk}t_{j,m}(Y_{kv}) = t_{j,m}(h_{v-u}^{(v-1)}(\mathbf{x}_m)).
\end{equation}
By Lemma \ref{lem_matrix_entries}, \eqref{eq_t_jm} is the identity
\begin{equation}
\label{eq_given}
\sum_k h_{k-u}^{(k-1)} (-1)^{k+v} t_{j,m}(e_{v-k}^{(k)}[\wh{m}]) = t_{j,m}(h_{v-u}^{(v-1)}(\mathbf{x}_m)).
\end{equation}
Setting $u = r+m-\alpha_i, v = r+m,$ and $k = r+s$, \eqref{eq_given} becomes
\begin{equation}
\label{eq_given2}
\sum_{s} h_{\alpha_i-m+s}^{(r+s-1)} (-1)^{s+m} t_{j,m}(e_{m-s}^{(r+s)}[\wh{m}]) = t_{j,m}(h_{\alpha_i}^{(r+m-1)}(\mathbf{x}_m)).
\end{equation}
Since $e_{m-s}^{(r+s)}[\wh{m}] = 0$ unless $0 \leq m-s \leq m-1$, we only need to sum over values of $s$ between $1$ and $m$. Thus, the two sides of \eqref{eq_given2} are precisely the $i,j$ entries of the two sides of \eqref{eq_ultimatewant}.

To complete the proof, observe that $H_{\delta}^{(r)}$ is upper triangular with 1's on the diagonal, so it has determinant 1, and thus \eqref{eq_ultimatewant} with $\alpha = \delta$ implies
\[
\det(E^{(r)}) = \det(A_{\delta}^{(r)}) = a_{\delta}^{(r)}.
\]
Taking determinants of \eqref{eq_ultimatewant}, we obtain
\[
s_{\lambda}^{(r)} a_{\delta}^{(r)} = a_{\lambda + \delta}^{(r)}.
\]
\end{proof}

\begin{remark}
If colors are identified (i.e., if $x_i^{(j)}$ is specialized to $y_i$ for all $j$) then the birational $S_m$ action reduces to the permutation of the variables $y_i$, so $(A_\alpha^{(r)})_{ij} = y_j^{\alpha_i}$, and Theorem \ref{thm_roa} specializes to the original definition of the Schur polynomials. Under this specialization, we also have $(B^{(r)}_{\alpha})_{ij} = y_j^{\alpha_i} \sigma_{(n-1)(m-j)}(y_j, \ldots, y_m)$, where $\sigma_i(y_j, \ldots, y_m)$ is the sum of monomials of degree $i$ in $y_j, \ldots, y_m$ such that the exponents of $y_{j+1}, \ldots, y_m$ are at most $n-1$. Thus, $b^{(r)}_{\alpha}$ becomes
\begin{equation}
\label{eq_Jucis_product}
a_{\alpha} \prod_{j=1}^m \sigma_{(n-1)(m-j)}(y_j, \ldots, y_m),
\end{equation}
and Theorem \ref{thm_roa_alt} is trivially equivalent to Theorem \ref{thm_roa} in this case.
\end{remark}

\begin{remark}
\label{rem_Jucis_formula}
There is some interesting combinatorics associated to the product in \eqref{eq_Jucis_product}. It is easy to see that
\[
\sigma_{(n-1)(m-j)}(y_j, \ldots, y_m) = \prod_{k=j+1}^m (y_j^{n-1} + y_j^{n-2}y_k + \cdots + y_j y_k^{n-2} + y_k^{n-1}),
\]
so
\begin{align*}
\prod_{j=1}^m \sigma_{(n-1)(m-j)}(y_j, \ldots, y_m) = \prod_{1 \leq j < k \leq m} \dfrac{y_j^n - y_k^n}{y_j - y_k} = \dfrac{a_{n\delta}}{a_{\delta}} = s_{(n-1)\delta}(y_1, \ldots, y_m).
\end{align*}
Jucis \cite{Jucis} gave a bijective proof of this identity using Schensted insertion. As mentioned in the proof of Theorem \ref{thm_roa_alt}, Lam and Pylyavskyy \cite{LPenergy} proved a ``colored refinement'' of this identity:
\[
\prod_{j=1}^m \sigma^{(r+j-1)}_{(n-1)(m-j)}(\mb{x}_j, \ldots, \mb{x}_m) = s^{(r)}_{(n-1)\delta}(\mb{x}_1, \ldots, \mb{x}_m).
\]
(They used this formula to show that the loop Schur function $s^{(n)}_{(n-1) \delta}$ tropicalizes to the intrinsic energy function for tensor products of one-row crystals of type $A_{n-1}^{(1)}$.) Their proof is algebraic; we would like to have a combinatorial proof, perhaps using a ``colored refinement'' of Schensted insertion. It would also be nice to have a combinatorial proof of Theorem \ref{thm_roa_alt}, perhaps along the lines of the argument in \cite{CRK}.
\end{remark}


\section{The loop Murnaghan--Nakayama rule}
\label{sec_mn}

The Murnaghan--Nakayama rule gives the Schur expansion of the product of a power sum symmetric function and a Schur function. A loop generalization of this rule was stated in \cite{ICCM}, and proved combinatorially by Ross \cite{Ross}. Here we give a short proof based on Theorem \ref{thm_roa}.

Recall the loop power sums $p_k = \sum_{i=1}^m \pi_i^k$, where $\pi_i = x_i^{(1)} x_i^{(2)} \cdots x_i^{(n)}$. Recall also that a \emph{ribbon} (or \emph{border strip}, or \emph{rim hook}) is a connected\footnote{Two boxes that share only a vertex are not considered to be connected; thus, the diagram $(2,1)/(1)$ is not connected.} skew Young diagram that does not contain any $2 \times 2$ squares. The \emph{size} of a ribbon is the number of boxes it contains, and the \emph{height} of a ribbon (denoted ${\rm ht}$) is one less than the number of rows it contains.

\begin{thm}
\label{thm_mn}
Let $\la$ be a partition with at most $m$ parts, and $k$ a positive integer. Then
\[
p_k(\mathbf{x}_1, \ldots, \mathbf{x}_m) s_{\lambda}^{(r)}(\mathbf{x}_1, \ldots, \mathbf{x}_m) = \sum (-1)^{{\rm ht}(\mu / \lambda)} s_{\mu}^{(r)}(\mathbf{x}_1, \ldots, \mathbf{x}_m),
\]
where the sum is over all partitions $\mu$ (with at most $m$ parts) such that $\mu / \lambda$ is a ribbon of size $kn$.
\end{thm}

\begin{proof}
All loop symmetric functions are in the variables $\mb{x}_1, \ldots, \mb{x}_m$. Since the birational $S_m$ action permutes the $\pi_i$, and each color appears the same number of times in $\pi_i$, we have
\begin{align*}
p_k a_{\lambda+\delta}^{(r)} &= \sum_{i=1}^m \pi_i^k \sum_{w \in S_m} {\rm sgn}(w) \prod_{j=1}^m t_{w(j),m}(x_m^{(r+m-1)} x_m^{(r+m-2)} \cdots x_m^{(r+m-(\lambda+\delta)_j)}) \\
&= \sum_{w \in S_m} {\rm sgn}(w) \sum_{i=1}^m \pi_{w(i)}^k \prod_{j=1}^m t_{w(j),m}(x_m^{(r+m-1)} x_m^{(r+m-2)} \cdots x_m^{(r+m-(\lambda+\delta)_j)}) \\
&= \sum_{w \in S_m} {\rm sgn}(w) \sum_{i=1}^m \prod_{j = 1}^m t_{w(j),m}(x_m^{(r+m-1)} x_m^{(r+m-2)} \cdots x_m^{(r+m-(\lambda+\delta)_j-kn\rho_{i,j})}) \\
&= \sum_{i=1}^m a_{\la + \delta + kn\epsilon_i}^{(r)},
\end{align*}
where $\rho_{i,j}$ is 1 if $i=j$ and 0 otherwise, and $\epsilon_i$ is the $i^{th}$ standard basis vector in $\mathbb{Z}^m$. Since $a_{\alpha}^{(r)}$ is anti-symmetric with respect to permuting parts of $\alpha$, a bit of bookkeeping (as in \cite[\S I.3, Ex. 11]{Mac}) shows that
\[
\sum_{i = 1}^m a_{\la + \delta + kn\epsilon_i}^{(r)} = \sum (-1)^{{\rm ht}(\mu/\lambda)} a_{\mu + \delta}^{(r)},
\]
where the sum is over all partitions $\mu$ (with at most $m$ parts) such that $\mu/\lambda$ is a ribbon of size $kn$.
Now divide by $a_\delta^{(r)}$ and apply Theorem \ref{thm_roa}.
\end{proof}

\bibliographystyle{alpha}
\bibliography{ROA_paper.refs}

\newcommand{\etalchar}[1]{$^{#1}$}
\begin{thebibliography}{HHI{\etalchar{+}}01}

\bibitem[BK00a]{BerKazh}
A.~Berenstein and D.~Kazhdan.
\newblock Geometric and unipotent crystals.
\newblock {\em Geom. Funct. Anal.}, Special Volume, Part I:188--236, 2000.
\newblock GAFA 2000 (Tel Aviv, 1999).

\bibitem[BK00b]{BravKazh}
A.~Braverman and D.~Kazhdan.
\newblock {$\gamma$}-functions of representations and lifting (with an appendix
  by {V}. {V}ologodsky).
\newblock {\em Geom. Funct. Anal.}, Special Volume, Part I:237--278, 2000.
\newblock With an appendix by V. Vologodsky, GAFA 2000 (Tel Aviv, 1999).

\bibitem[CRK95]{CRK}
J.~Carbonara, J.~Remmel, and A.~Kulikauskas.
\newblock A combinatorial proof of the equivalence of the classical and
  combinatorial definitions of {S}chur function.
\newblock {\em J. Combin. Theory Ser. A}, 72(2):293--301, 1995.

\bibitem[Eti03]{Et}
P.~Etingof.
\newblock Geometric crystals and set-theoretical solutions to the quantum
  {Y}ang--{B}axter equation.
\newblock {\em Comm. Algebra}, 31(4):1961--1973, 2003.

\bibitem[Fri]{FriR}
G.~Frieden.
\newblock The geometric {$R$}-matrix for affine crystals of type {$A$}.
\newblock Preprint: arXiv 1710.07243.

\bibitem[HHI{\etalchar{+}}01]{HHIKTT}
G.~Hatayama, K.~Hikami, R.~Inoue, A.~Kuniba, T.~Takagi, and T.~Tokihiro.
\newblock The {$A_M^{(1)}$} automata related to crystals of symmetric tensors.
\newblock {\em J. Math. Phys.}, 42(1):274--308, 2001.

\bibitem[Juc80]{Jucis}
A.-A.~A. Jucis.
\newblock Tournaments and generalized {Y}oung tableaux.
\newblock {\em Mat. Zametki}, 27(3):353--359, 492, 1980.
\newblock (English translation: Math. Notes 27(3-4):175--178, 1980).

\bibitem[KNO10]{KNO}
M.~Kashiwara, T.~Nakashima, and M.~Okado.
\newblock Tropical {$R$} maps and affine geometric crystals.
\newblock {\em Represent. Theory}, 14:446--509, 2010.

\bibitem[Lam12]{ICCM}
T.~Lam.
\newblock Loop symmetric functions and factorizing matrix polynomials.
\newblock In {\em {F}ifth {I}nternational {C}ongress of {C}hinese
  {M}athematicians. {P}art 1, 2}, volume~2 of {\em {AMS/IP} Stud. Adv. Math.,
  51, pt. 1}, pages 609--627. Amer. Math. Soc., Providence, RI, 2012.

\bibitem[LP]{LPplus}
T.~Lam and P.~Pylyavskyy.
\newblock Private communication.

\bibitem[LP12]{LPwhirl}
T.~Lam and P.~Pylyavskyy.
\newblock Total positivity in loop groups, {I}: {W}hirls and curls.
\newblock {\em Adv. Math.}, 230(3):1222--1271, 2012.

\bibitem[LP13a]{LPnet}
T.~Lam and P.~Pylyavskyy.
\newblock Crystals and total positivity on orientable surfaces.
\newblock {\em Selecta Math. (N.S.)}, 19(1):173--235, 2013.

\bibitem[LP13b]{LPenergy}
T.~Lam and P.~Pylyavskyy.
\newblock Intrinsic energy is a loop {S}chur function.
\newblock {\em J. Comb.}, 4(4):387--401, 2013.

\bibitem[LPS]{LPS}
T.~Lam, P.~Pylyavskyy, and R.~Sakamoto.
\newblock Rigged configurations and cylindric loop {S}chur functions.
\newblock {\em Annales de l'Institut Henri Poincar\'e D}, to appear.
\newblock Preprint: arXiv 1410.4455.

\bibitem[Mac95]{Mac}
I.~G. Macdonald.
\newblock {\em Symmetric functions and {H}all polynomials}.
\newblock Oxford Mathematical Monographs. The Clarendon Press, Oxford
  University Press, New York, second edition, 1995.
\newblock With contributions by A. Zelevinsky, Oxford Science Publications.

\bibitem[Ros14]{Ross}
D.~Ross.
\newblock The loop {M}urnaghan--{N}akayama rule.
\newblock {\em J. Algebraic Combin.}, 39(1):3--15, 2014.

\bibitem[RZ13]{RossZong}
D.~Ross and Z.~Zong.
\newblock The gerby {G}opakumar--{M}ari\~no--{V}afa formula.
\newblock {\em Geom. Topol.}, 17(5):2935--2976, 2013.

\bibitem[Sta12]{EC1}
R.~P. Stanley.
\newblock {\em Enumerative combinatorics. {V}olume 1}, volume~49 of {\em
  Cambridge Studies in Advanced Mathematics}.
\newblock Cambridge University Press, Cambridge, second edition, 2012.

\bibitem[Yam01]{Yam}
Y.~Yamada.
\newblock A birational representation of {W}eyl group, combinatorial
  {$R$}-matrix and discrete {T}oda equation.
\newblock In {\em Physics and combinatorics, 2000 (Nagoya)}, pages 305--319.
  World Sci. Publ., River Edge, NJ, 2001.

\end{thebibliography}

\end{document}